\documentclass[a4paper,twoside,12pt,reqno]{amsart}
\usepackage{amssymb,amsfonts,amscd,euscript,latexsym}
\usepackage{xy} \xyoption{all}

\newtheorem{lemma}{Lemma}[section]
\newtheorem{proposition}[lemma]{Proposition}
\newtheorem{theorem}[lemma]{Theorem}
\newtheorem{corollary}[lemma]{Corollary}

\theoremstyle{remark}
\newtheorem{remark}[lemma]{Remark}

\theoremstyle{plain}
\newenvironment{proofof[1]}{\em Proof of {#1}.\quad}{\quad\qed}

\numberwithin{equation}{section}

\newcommand{\NN}{{\mathbb N}}

\newcommand{\eps}{{\varepsilon}}

\newcommand{\ol}{\overline}
\newcommand{\mloc}{{M_{\text{\rm loc}}(A)}}

\newcommand{\Mloc}[1]{{M_{\text{\rm loc}}({#1})}}

\newcommand{\qmax}{{Q_{\text{\rm max}}(A)}}
\newcommand{\Qmax}[1]{{Q_{\text{\rm max}}({#1})}}

\newcommand{\Qmaxs}[1]{{Q_{\text{\rm max}}^s({#1})}}

\newcommand{\Qmaxsb}[1]{{Q_{\text{\rm max}}^s({#1})_b}}
\newcommand{\qmaxsb}{{\Qmaxsb A}}

\newcommand{\CBsAA}{{C\!B_{A\text{-}A}}}

\newcommand{\longrightarrowraised}{{\hbox{\raise.5\jot\hbox{\scriptsize$\longrightarrow$}}}}
\newcommand{\longleftarrowraised}{{\hbox{\raise.5\jot\hbox{\scriptsize$\longleftarrow$}}}}
\newcommand{\dirlim}{{\smash{\underset{\longrightarrowraised}
                   {\operatorname{lim}}}\vphantom{a^{}_{a_f^{}}}}}

\newcommand{\Dirlim}[1]{{\smash{\dirlim{}}\sp{}_{\,{#1}}
                   \vphantom{a^{}_{a_f^{}}}}}
\newcommand\alglim{{\smash{\underset{\longrightarrowraised}
                   {\operatorname{alg\,lim}}}\vphantom{a^{}_{a_f^{}}}}}
\newcommand\Alglim[1]{{\smash{\alglim{}}\sp{}_{\,{#1}}
                   \vphantom{a^{}_{a_f^{}}}}}

\newcommand{\Oone}{{\mathcal O}_1}

\newcommand{\Icer}{{{\mathfrak I}_{cer}}}

\newcommand{\id}[1]{{\text{\rm id}_{{#1}}}}
\newcommand{\Abar}{{\hbox{\kern.22em$\ol{\phantom{I}}\kern-.75em A$\kern.09em}}}
\newcommand{\smallAbar}{{\hbox{\scriptsize\kern.22em$\ol{\phantom{I}}\kern-.75em A$\kern.09em}}}

\newcommand{\Bbar}{{\hbox{\kern.22em$\ol{\phantom{I}}\kern-.75em B$\kern.09em}}}

\newcommand{\cbnorm}[1]{{\|{#1}\|_{cb}}}

\newcommand{\rest}[2]{{{#1}_{\kern-.5pt|{#2}}}}  
\newcommand{\haaga}{{A\otimes_h\mkern-2mu A}}
\def\C*{{\sl C*}-algebra}
\def\Cs*{{\sl C*}-subalgebra}
\def\W*{{\sl W*}-algebra}
\def\AW*{{\sl AW*}-algebra}
\def\AF/{{\sl AF}-algebra}
\def\AFs/{{\sl AF}-algebras}
\def\vN/{{von Neumann algebra}}
\def\sigfinite/{{$\sigma$-finite}}

\begin{document}

\title[The Maximal C*-Algebra of Quotients as an Operator Bimodule]%
      {The Maximal C*-Algebra of Quotients\\ as an Operator Bimodule}
\author{Pere Ara, Martin Mathieu and Eduard Ortega}
\address{Departament de Matem\`atiques, Universitat Aut\`onoma de Barcelona,
         E-08193 Bellaterra (Barcelona), Spain ({\rm P. Ara})}
\email{para@mat.uab.cat}
\address{Department of Pure Mathematics, Queen's University Belfast, Belfast BT7 1NN,
         Northern Ireland ({\rm M. Mathieu})}
\email{m.m@qub.ac.uk}
\address{Department of Mathematics and Computer Science, University of Southern Denmark,
         Campusvej~55, DK-5230 Odense~M, Denmark ({\rm E. Ortega})}
\email{ortega@imada.sdu.dk}
\thanks{The first-named and the third-named author's research was partially supported by the DGI and the European Regional Fund, jointly,
through project MTM2005-00934 and, in addition, by the Comissionat per Universitats i Recerca de la Generalitat de Catalunya.
The second-named author was partially supported by a Scheme~4 grant of the London Mathematical Society.
Part of this paper was written during mutual visits of the second-named author to the University of Southern
Denmark, Odense and the third-named author to Queen's University Belfast and both would like to thank the
respective Mathematics Departments for their hospitality.}


\begin{abstract}
We establish a description of the maximal \C* of quotients of a unital \C*~$A$
as a direct limit of spaces of completely bounded bimodule homomorphisms from certain
operator submodules of the Haagerup tensor product $\haaga$ labelled by the essential
closed right ideals of $A$ into~$A$. In addition the invariance of the construction of the
maximal \C* of quotients under strong Morita equivalence is proved.
\end{abstract}

\maketitle

\section{Introduction}\label{sec:intro}
The maximal \C* of quotients, $\qmax$, of a (unital) \C* $A$ was introduced in~\cite{AM08}
as a {\sl C*}-analytic analogue of the maximal symmetric ring of quotients of a non-singular ring
studied, e.g., in~\cite{Lann} and~\cite{Ed06}. As a \C* of quotients it shares some of the properties
of the local multiplier algebra $\mloc$ of~$A$~\cite{AMMzth}; for instance, it arises as
the completion of the bounded part of its algebraic counterpart and it can be canonically embedded into the injective envelope
$I(A)$ of~$A$. However, in contrast to the situation of $\mloc$, there is no direct limit construction for
$\qmax$ in the category of \C*s. This makes its study somewhat more cumbersome. The purpose of this
note is to alleviate this difficulty by providing a direct limit description in a different category,
the category of operator modules. From the well-established good properties of the Haagerup tensor
product it comes as no surprise that this concept will play an important role.

We will be guided by the construction of the maximal symmetric ring of quotients $\Qmaxs R$ of a (non-singular unital) ring $R$ with
involution as a two-sided localisation of the regular bimodule ${}_RR_R$ with respect to a certain filter of right ideals in~\cite{Ed08a}.
In this situation, for each essential right ideal $I$ of $R$, a certain right ideal $M_I$ in the ring $R\otimes R^{\kern.2\jot\rm op}$ is used to
introduce the filter $\Omega$ consisting of all right ideals containing some~$M_I$. It then turns out that $\Qmaxs R$ is canonically
isomorphic to $\Dirlim\Omega{\rm Hom}(M_I,R)$, where ${\rm Hom}(M_I,R)$ stands for the space of all $R$-bimodule homomorphisms
from $M_I$ into~$R$. Our principal goal in the present paper is to modify this algebraic construction in such a way that, in the situation
of a unital \C*~$A$, the \C* $\qmax$ can be obtained in an analogous manner as a direct limit of spaces of completely bounded bimodule homomorphisms
from certain operator submodules of the Haagerup tensor product $\haaga$ labelled by the essential closed right ideals of $A$ into~$A$.

On our route to establish our main result, Theorem~\ref{thm:qmax-direct-limit}, we shall obtain several other useful properties of the maximal  \C*
of quotients. These will enable us to show, in Theorem~\ref{thm:morita-invariance}, that this construction is invariant under
strong Morita equivalence of \C*s (as it is indeed the case for the local multiplier algebra).  For a comprehensive
discussion of various types of Morita equivalence, we refer to~\cite{BMP}.

Throughout, we shall use the terminology and notation of~\cite{AMMzth}, \cite{AM08}
and~\cite{ER} with the exception that $\haaga$ already denotes the completion of the algebraic tensor product $A\otimes A$ with
respect to the Haagerup norm~$\|\cdot\|_h$. In particular, $\Oone$ stands for the category of operator spaces with complete contractions
as the morphisms, and $\CBsAA(E,F)$ is the operator space of all completely bounded $A$-bimodule maps from the operator
$A$-bimodule $E$ into the operator $A$-bimodule~$F$.

\section{The Results}\label{sec:results}

For our purposes here, the following description of the maximal \C* of quotients as a \Cs* of the injective envelope $I(A)$
of a \C* $A$ is the most expedient. Let $\qmaxsb=\{q\in I(A)\mid qJ+q^*J\subseteq A\text{\ for some }J\in\Icer\}$, where $\Icer$
denotes the set of all closed essential right ideals of~$A$. Then $\qmax=\ol\qmaxsb$; see \cite[Theorem~4.8]{AM08}.

The proposition below states a property of the maximal \C* of quotients that is
shared with the maximal symmetric ring of quotients of a semiprime ring, compare \cite[Proposition~1.6]{Ed06}.

\begin{lemma}\label{lem:ideals-in-corner}
Let\/ $A$ be a \C*, and let\/ $e\in M(A)$ be a projection. Then
\[
\iota \colon\Icer(eAe)\longrightarrow \Icer(A),\quad J\longmapsto \ol{JA+(1-e)A}
\]
is an injective mapping and
\[
\kappa \colon \Icer(A)\longrightarrow\Icer(eAe),\quad I\longmapsto eIe
\]
is a surjective mapping such that\/ $\kappa\circ\iota =\id{\Icer(eAe)}$.
\end{lemma}
\begin{proof} 
Take $J\in \Icer(eAe)$ and note that the annihilator $(AeA)^{\perp}$ coincides with $(\ol{AeA})^{\perp}$ and is contained
in $(1-e)A(1-e)$. In order to prove that $\ol{JA+(1-e)A}$ is an essential right ideal of $A$, it thus suffices to show that
$\ol{JA+(1-e)AeA}+(AeA)^{\perp}$ is an essential right ideal of~$A$. Since $\ol{AeA}\oplus(AeA)^{\perp}$ is essential as a right
ideal, it is enough to show that $\ol{JA+(1-e)AeA}$ is essential
in $\ol{AeA}$. Let $z$ be a non-zero element in $\ol{AeA}$. If $ezAe=0$ then $ez\ol{AeA}=0$ and thus $ez=0$ entailing that
$z=(1-e)z\in \ol{(1-e)AeA}$. We can therefore assume that there is
$t\in A$ such that $ezte$ is a non-zero element of $eAe$. Since $J$
is essential in $eAe$, there is $s\in eAe$ such that $eztes$ is a non-zero element of~$J$. 
It follows that
\[
z(tes)=eztes+(1-e)ztes\in J+\ol{(1-e)AeA}\subseteq \ol{JA+(1-e)AeA}\, ,
\]
as desired. 

To prove the second assertion, take $I\in \Icer(A)$. If $z$ is a non-zero element in $eAe$, there is
$a\in A$ such that $za$ is a non-zero element of~$I$. Since $za\in AeA$, we have $zaAe\ne 0$; thus 
$0\ne z(eaAe)\subseteq eIe$ wherefore $eIe\in \Icer(eAe)$. It is obvious
that $\kappa\circ \iota =\id{\Icer(eAe)}$. In particular, $\iota$ is injective and $\kappa$ is surjective.
\end{proof}

\begin{proposition}\label{prop:corner-of-qmax}
Let\/ $e$ be a projection in the multiplier algebra of a \C*~$A$. Then\/ $\Qmax{eAe}=e\qmax e$.
\end{proposition}
\begin{proof} \quad
By  \cite[Proposition 6.3]{Hametb}, $I(eM(A)e)=eI(M(A))e$; combining this with
$I(M(A))=I(A)$ and $I(M(eAe))=I(eAe)$, we get
\[
I(eAe)=I(M(eAe))=I(eM(A)e)=eI(M(A))e=eI(A)e.
\]
Thus we can view $\Qmax{eAe}$ as a \Cs* of $eI(A)e$. We shall divide the proof of the main statement into three steps.

\smallskip\noindent
Step~1. \qquad $e\qmaxsb e\subseteq \Qmaxsb{eAe}$.

\smallskip
In order to show this take $q\in\qmaxsb$. There is $I\in\Icer(A)$
such that $eqeI$ and $eq^*eI$ are both contained in~$A$. By
Lemma~\ref{lem:ideals-in-corner}, $eIe\in\Icer(eAe)$ and, since
$eqe\,eIe\subseteq eAe$ and $eq^*e\,eIe\subseteq eAe$, and $eqe\in
I(eAe)$, we conclude that $eqe\in \Qmaxsb{eAe}$.

\smallskip\noindent
Step~2.\qquad$\Qmaxsb{eAe}\subseteq e\qmaxsb e$.

\smallskip
To see this, let $q\in\Qmaxsb{eAe}\subseteq I(eAe)=eI(A)e\subseteq
I(A)$ and take $J\in\Icer(eAe)$ such that $qJ$, $q^*J\subseteq eAe$.
Putting $J'=\ol{JA+(1-e)A}$ we obtain an essential right ideal
of~$A$ (Lemma~\ref{lem:ideals-in-corner}). Note that
$$qJ'+q^*J'\subseteq \ol{qJA}+\ol{q^*JA}\subseteq A$$
and thus $q\in e\qmaxsb e$.

\smallskip\noindent
Step~3.\quad From the first two steps we conclude, taking closures
in $I(A)$, that
\begin{equation*}
e\qmax e = e\ol{\qmaxsb}e = \ol{e\qmaxsb e} = \ol{\Qmaxsb{eAe}} = \Qmax{eAe},
\end{equation*}
since $e\qmax e$ is closed.
\end{proof}

\begin{corollary}\label{cor:matrices-over-qmax}
Let\/ $A$ be a unital \C*. Then\/ $M_n(\qmax)=\Qmax{M_n(A)}$ for each\/ $n\in\NN$.
\end{corollary}
\begin{proof}
This follows from a standard argument, see, e.g., \cite[Remark~17.6]{Lam98}, so we merely sketch the essential part for completeness.

Let $\{e_{ij}\mid1\leq i,j\leq n\}$ be the canonical matrix units in $M_n(A)\subseteq\Qmax{M_n(A)}$ and denote by $e_i=e_{ii}$ the
$n$ mutually orthogonal projections with $\sum_{i=1}^ne_i=1$. For $x\in\Qmax{M_n(A)}$ set
$a_{ij}=\sum_{m=1}^n e_{mi}\,x\,e_{jm}$. Then $x=\sum_{i,j=1}^na_{ij}\,e_{ij}$, since
\begin{equation}\label{eq:commuting}
a_{ij}\,e_{k\ell}=e_{ki}\,x\,e_{j\ell}=e_{k\ell}\,a_{ij}\qquad(1\leq i,j,k,\ell\leq n)
\end{equation}
and thus $a_{ij}\,e_{ij}=e_{ii}\,x\,e_{jj}=e_ixe_j$. Letting $q_{ij}^{(m)}=e_{mi}\,x\,e_{jm}$ we infer from~(\ref{eq:commuting})
that $q_{ij}^{(m)}=q_{ij}$ for all $1\leq m\leq n$. By Proposition~\ref{prop:corner-of-qmax},
\[
q_{ij}=e_m\,q_{ij}\,e_m\in e_m\,\Qmax{M_n(A)}\,e_m=\Qmax{e_m\,M_n(A)\,e_m}=\qmax
\]
as each $e_m$ is a full projection in $M_n(A)$. Thus, with $q_x=\sum_{i,j}q_{ij}\,e_{ij}\in M_n(\qmax)$, we establish a *-isomorphism
$x\mapsto q_x$ from $\Qmax{M_n(A)}$ onto $M_n(\qmax)$ (which also satisfies $\Qmaxsb{M_n(A)}=M_n(\qmaxsb)$).
\end{proof}

We obtain the following very useful consequence.

\begin{theorem}\label{thm:morita-invariance}
Let\/ $A$ and\/ $B$ be two unital strongly Morita equivalent \C*s. Then their maximal \C*s of quotients\/ $\qmax$ and\/
$\Qmax B$ are strongly Morita equivalent.
\end{theorem}
\begin{proof}
Suppose that $A$ and $B$ are strongly Morita equivalent. Then there is a full projection $e$ in some matrix algebra $M_n(A)$
such that $B\cong eM_n(A)e$. Therefore
\begin{equation}\label{eq:for-morita}
\Qmax{B}\cong\Qmax{eM_n(A)e}=e\,\Qmax{M_n(A)}\,e=e\,M_n(\qmax)\,e
\end{equation}
by Proposition~\ref{prop:corner-of-qmax} and Corollary~\ref{cor:matrices-over-qmax}, since $e$ is a full projection in
$\Qmax{M_n(A)}$. As a result, $\Qmax B$ and $\qmax$ are strongly Morita equivalent.
\end{proof}

\begin{corollary}\label{cor:stability}
Let\/ $A$ and\/ $B$ be two unital strongly Morita equivalent \C*s. If\/ $A=\qmax$ then\/ $B=\Qmax B$.
\end{corollary}

\begin{remark}\label{rem:morita-for-mloc}
It is evident that the analogues of Theorem~\ref{thm:morita-invariance} and its corollary hold for the local multiplier algebra
once the analogue of the statement in Proposition~\ref{prop:corner-of-qmax} is verified. In fact,
$\Mloc{eAe}=e\mloc e$ for any projection $e$ in a unital \C*~$A$. If $e$ is full, then $A$ and the hereditary \Cs* $eAe$ are
strongly Morita equivalent and therefore the lattices of their closed essential ideals are isomorphic \cite[Proposition~1.2.38]{AMMzth}.
Since $eM(I)e=M(eIe)$ for every closed essential ideal $I$ \cite[Corollary~1.2.37]{AMMzth}, the equality
$\Mloc{eAe}=e\mloc e$ follows from the direct limit formula for the local multiplier algebra. If $e$ is arbitrary, then, using
\cite[Proposition~2.3.6 ]{AMMzth}, we have
\begin{equation*}
e\mloc e=e\Mloc{J\oplus J^\perp}e=e\Mloc J e\oplus e\Mloc{J^\perp}e=e\Mloc J e=\Mloc{eAe},
\end{equation*}
where $J=\ol{AeA}$ is the closed ideal generated by $e$, $J^\perp$ is its annihilator (so that $J+J^\perp$ is essential) and we apply the above
argument to the full hereditary \Cs* $eAe$ of~$J$.

The Morita invariance of the local multiplier algebra has several pleasant consequences. For instance, any unital \C* $B$ which is strongly
Morita equivalent to a boundedly centrally closed unital \C* $A$ (i.e., $Z(A)=Z(\mloc)$ \cite[Definition~3.2.1]{AMMzth})
is boundedly centrally closed. For, the centres of $A$ and $B$ and the centres of $\mloc$ and $\Mloc B$ are isomorphic.
In addition, the analogue of Corollary~\ref{cor:stability} can be used to study iterated local multiplier algebras. For example,
suppose that $A$ is strongly Morita equivalent to a commutative \C*~$B$. Since
$\Mloc{\Mloc B}=\Mloc B$ as $\Mloc B$ is an \AW*, it follows that $\Mloc{\mloc}=\mloc$.

We like to point out that the above arguments are restricted to the situation of \textit{unital\/} \C*s, as the example
$A=C[0,1]\otimes K(\ell^2)$ studied in \cite{AM08} and \cite{ArgFar2} shows; there, $\Mloc\mloc\neq\mloc$.
\end{remark}

We now turn our attention to the description of the maximal \C* of quotients as a direct limit of spaces of completely
bounded bimodule homomorphisms. To this end
we shall consider the Haagerup tensor product $\haaga$ of a unital \C* $A$ as an operator $A$-bimodule with the operations
$a(x\otimes y)b=ax\otimes yb$, $a,b\in A$ and $x,y\in A$.
For an essential closed right ideal $I$ in $A$, we define the closed sub-bimodule
\[
M_I=\overline{A\otimes I+I^*\otimes A}.
\]
If $I$ is two-sided, the closure is not needed \cite[Theorem~3.8]{ASS93} but we do not have
an analogous result available for one-sided ideals.

\begin{proposition}\label{prop:cbmaps_on_MI}
Let\/ $A$ be a unital \C* and let\/ $I\in\Icer$. For each completely bounded
$A$-bimodule homomorphism\/ $\psi\colon M_I\to A$ there exists a unique\/ $q\in\qmaxsb$
such that\/ $\rest{\psi}{1\otimes I}=\rest{L_q}{I}$ and\/ $\cbnorm\psi=\|q\|$.
\end{proposition}

We prepare the proof by the following simple lemma.

\begin{lemma}\label{lem:cbmaps on aha}
Let\/ $A$ be a unital \Cs* of the unital \C*\/ $B$. For each\/ $q\in B$, by\/
$\psi_q(a\otimes b)=aqb$, $a,b\in A$ we can define a completely bounded $A$-bimodule homomorphism\/
$\psi_q\colon\haaga\to B$ such that\/ $\cbnorm{\psi_q}=\|q\|$.
\end{lemma}
\begin{proof}
Evidently, the assignment $a\otimes b\mapsto aqb$ yields an $A$-bimodule homomorphism $\psi_q$ from $A\otimes A$ into~$B$.
Since $\psi_q(1\otimes1)=q$, we have $\|\psi_q\|\geq\|q\|$. For the reverse inequality, take
$u\in M_n(A\otimes A)$ with $\|u\|_h=1$. Let $\eps>0$.
Then there exist $v\in M_{n,r}(A)$, $w\in M_{r,n}(A)$, for some $r\in\NN$, such that $u=v\odot w$ and $\|v\|\,\|w\|<1+\eps$.
If $v=(v_{k\ell})$, $w=(w_{k\ell})$ then
\begin{equation*}
\psi_q^{(n)}(u)= \psi_q^{(n)}(v\odot w)
                               =  \sum_{\ell=1}^r\bigl(\psi_q(v_{i\ell}\otimes w_{\ell j})\bigr)_{ij}
                               = \sum_{\ell=1}^r\bigl(v_{i\ell}\,q\,w_{\ell j}\bigr)_{ij}
                               =(v_{k\ell})\cdot q\cdot(w_{k\ell}),
\end{equation*}
where, by abuse of notation, we also denote by $q$ the $r\times r$ diagonal matrix with $q$ along the diagonal.
Therefore,
$\|\psi_q^{(n)}(u)\|\leq\|v\|\,\|q\|\,\|w\|<(1+\eps)\,\|q\|$. It follows that $\|\psi_q^{(n)}\|\leq\|q\|$ for all $n\in\NN$ wherefore
$\cbnorm{\psi_q}=\|q\|$. Extending $\psi_q$ to $\haaga$ thus yields the result.
\end{proof}

\noindent
\begin{proofof}[Proposition~\ref{prop:cbmaps_on_MI}]\rm
Let $x\in I$ and define $g(x)=\psi(1\otimes x)$ and
$f(x^*)=\psi(x^*\otimes 1)$. Then $g\colon I\to A$ and $f\colon I^*\to A$ are completely
bounded right and left module homomorphisms, respectively, such that
$y^*g(x)=f(y^*)x$ for all $x,y\in I$. By \cite[Lemma~4.7]{AM08} and the subsequent remarks, there is a unique element $q\in I(A)$ such that
$g=\rest{\psi}{1\otimes I}=\rest{L_q}{I}$ and $f=\rest{\psi}{I^*\otimes 1}=\rest{R_q}{I^*}$. As $qI+q^*I\subseteq A$, it follows
that $q\in\qmaxsb$ and $\|q\|=\|{L_q}_{|I}\|\leq \cbnorm\psi$.

For $u=a\otimes x+y^*\otimes b\in M_I$, $x,y\in I$, $a,b\in A$ we thus obtain
\[
\psi(u)=\psi(a\otimes x+y^*\otimes b)=aqx+y^*qb.
\]

Apply Lemma~\ref{lem:cbmaps on aha} with $B=\qmax$ to obtain $\psi_q$ on $\haaga$ with the property that
$\rest{\psi_q}{M_I}=\psi$. It follows that $\|q\|\leq\cbnorm\psi\leq\cbnorm{\psi_q}=\|q\|$, as desired.
\end{proofof}

\medskip

For $I,J\in\Icer$ with $J\subseteq I$ we have an embedding $M_J\subseteq M_I$ and thus we can define
\begin{equation*}
\rho_{JI}\colon\CBsAA(M_I,A)\longrightarrow\CBsAA(M_J,A),\quad\psi\longmapsto\rest{\psi}{M_J}.
\end{equation*}
The fundamental property that $\|\rest{L_q}{J}\|=\|\rest{L_q}{I}\|=\|q\|$ together with
Proposition~\ref{prop:cbmaps_on_MI} entails that
$\cbnorm\psi=\|\psi_{|1\otimes I}\|=\|\psi_{|1\otimes J}\|=\|q\|$; i.e., each $\rho_{JI}$
is isometric. To show that $\rho_{JI}$ is indeed completely isometric, that is,
\[
\rho^{(n)}_{JI}\colon M_n(\CBsAA(M_I,A))\longrightarrow M_n(\CBsAA(M_J,A)),\quad
(\psi_{k\ell})\longmapsto(\rest{\psi_{k\ell}}{M_J}),
\]
$I,J\in\Icer$, $J\subseteq I$ is isometric for each $n\in\NN$,
we use that $M_n(\CBsAA(M_I,A))=\CBsAA(M_I,M_n(A))$
and that every completely bounded $A$-bimodule map $\psi\colon M_I\to M_n(A)$
gives rise to a unique element $q\in M_n(\qmaxsb)$ via $q=(q_{k\ell})$, where
$q_{k\ell}$ is the element in $\qmaxsb$ with $\rest{\psi_{k\ell}}{1\otimes I}=\rest{L_{q_{k\ell}}}{I}$ given by
Proposition~\ref{prop:cbmaps_on_MI}. Letting $\psi_{q_{k\ell}}$ denote the extension of $\psi_{k\ell}$ to
$\haaga$ with values in $\qmax$ as above we take $u\in A\otimes A$ with $\|u\|_h=1$. Given $\eps>0$, write it as $u=v\odot w$
for some $v\in M_{1,r}(A)$, $w\in M_{r,1}(A)$ such that $\|v\|\,\|w\|<1+\eps$.
Denoting by $q_{k\ell}$ once again the $r\times r$ diagonal matrix with $q_{k\ell}$ along the diagonal we find
\begin{equation*}
\bigl(\psi_{q_{k\ell}}(u)\bigr) =\bigl(v\cdot q_{k\ell}\cdot w\bigr)
                                                   =\begin{pmatrix}v&0&\cdots&0\\
                                                                                      0&\ddots&&\vdots\\
                                                                                     \vdots & &\ddots&0\\
                                                                                     0&\cdots&0&v
                                                     \end{pmatrix}
                                                     \begin{pmatrix}{q}_{11}&\cdots&{q}_{1n}\\
                                                                                      \vdots&&\vdots\\
                                                                                      {q}_{n1}&\cdots&{q}_{nn}
                                                     \end{pmatrix}
                                                     \begin{pmatrix}w&0&\cdots&0\\
                                                                                     0  &\ddots&&\vdots\\
                                                                                     \vdots & &\ddots&0\\
                                                                                    0&\cdots&0&w
                                                     \end{pmatrix}\\
\end{equation*}
entailing
\[
\bigl\|\bigl(\psi_{q_{k\ell}}(u)\bigr)\bigr\|
\leq\bigl\|v\bigr\|\,\bigl\|w\bigr\|\,\left\| \begin{pmatrix}(q_{k\ell})&0&\cdots&0\\
                                                                                                                        0               &\ddots&&\vdots\\
                                                                                                                         \vdots     & &\ddots&0\\
                                                                                                                        0&\cdots&0&(q_{k\ell})
                                                                                                                        \end{pmatrix} \right\|,
\]
by the canonical shuffle.

As before, this implies that the norm of the mapping $\psi=(\psi_{k\ell})$ is dominated by $\|(q_{k\ell})\|$, wherefore, in
$M_n(\CBsAA(M_I,A))$, $\|(\psi_{k\ell})\|\leq\|(q_{k\ell})\|$. (Note that, by Corollary~\ref{cor:matrices-over-qmax}, the mapping
$\psi_q$ given by Lemma~\ref{lem:cbmaps on aha} applied to $B=\Qmax{M_n(A)}$ agrees with $(\psi_{q_{k\ell}})$.)

On the other hand, from the fact that $\rest{\psi_{k\ell}}{1\otimes I}=\rest{L_{q_{k\ell}}}{I}$, we conclude that
\begin{equation*}
\|(\psi_{k\ell})\|\geq\|(\rest{\psi_{k\ell}}{1\otimes I})\|=\|(\rest{L_{q_{k\ell}}}{I})\|=\|(q_{k\ell})\|
\end{equation*}
by Lemma~3.9 in~\cite{AM08}.
As a consequence, the norm of $\psi=(\psi_{k\ell})\in M_n(\CBsAA(M_I,A))$ coincides with the norm of $q=(q_{k\ell})\in M_n(\qmaxsb)$
which, in particular, implies that the restriction homomorphisms $\rho^{(n)}_{JI}$ are isometric for each~$n$.

We are now in a position to prove our main result.

\begin{theorem}\label{thm:qmax-direct-limit}
Let\/ $A$ be a unital \C*. Then the operator $A$-bimodules\/ $\qmaxsb$ and\/ $\Alglim{I\in\Icer}\CBsAA(M_I,A)$ are completely
isometrically isomorphic. As a result, $\qmax=\Dirlim{I\in\Icer}\CBsAA(M_I,A)$ in the category\/~$\Oone$.
\end{theorem}
\begin{proof}
By the above arguments, the mapping
\[
\tau_I\colon\CBsAA(M_I,A)\longrightarrow\qmaxsb,\quad\psi\longmapsto q_\psi,
\]
where $q_\psi\in\qmaxsb$ is the unique element determined by Proposition~\ref{prop:cbmaps_on_MI}, is completely isometric
for each $I\in\Icer$. From the construction, it is clear that, for $J\subseteq I$, the following diagram is commutative
\begin{equation*}
\xymatrix{\CBsAA(M_I,A)\ar[rd]^{\tau_I}\ar[rr]^{\rho_{JI}} & & \CBsAA(M_J,A)\ar[ld]^{\tau_J} \\
 & \qmaxsb&  }
\end{equation*}
Since each $\rho_{JI}$ is completely isometric, there exists a complete isometry
\[
\tau\colon\Alglim{I\in\Icer}\CBsAA(M_I,A)\to\qmaxsb
\]
such that $\rest{\tau}{\CBsAA(M_I,A)}=\tau_I$. Let $q\in\qmaxsb$ and choose $I\in\Icer$ with the property that
$qI$, $q^*I\subseteq A$. Then
\[
\sum^r_{i=1}a_i\otimes x_i+\sum^s_{j=1}y^*_j\otimes b_j\longmapsto\sum^r_{i=1} a_iqx_i+\sum^s_{j=1} y^*_jqb_j
\]
defines a completely bounded $A$-bimodule homomorphism $\psi_q$ which clearly satisfies $\tau_I(\psi_q)=q$.
This yields a right inverse of $\tau$, which is therefore surjective. This complete isometry extends to a complete
isometry from $\Dirlim{I\in\Icer}\CBsAA(M_I,A)$ onto $\qmax$.

Although each individual space $\CBsAA(M_I,A)$ is not necessarily an $A$-bimodule, the direct limit
$\Alglim{I\in\Icer}\CBsAA(M_I,A)$ is. To see this note that, whenever $\psi\in\CBsAA(M_I,A)$ is given as $\psi=\psi_q$ with
$qI$, $q^*I\subseteq A$ and $c\in A$, we can define $\psi c\in\CBsAA(M_J,A)$ via
\begin{equation*}
(\psi c) (a\otimes x+y^*\otimes b)=aqcx+y^*qcb\qquad(a,b\in A,\,x,y\in J),
\end{equation*}
where $J\in\Icer$ is given by $J=\{u\in I\mid cu\in I\}$. A similar expression yields $c\psi$, and these module operations are
evidently compatible with the connecting maps $\rho_{JI}$. Therefore, we obtain an $A$-bimodule structure on
$\Alglim{I\in\Icer}\CBsAA(M_I,A)$ turning it into an operator $A$-bimodule, which is completely isometrically
isomorphic to $\qmaxsb$.
\end{proof}

\begin{remark}
What is the involution on $\qmax$ in the above picture? Endow $\haaga$ with an involution ${}^*$ defined by
$(x\otimes y)^*=y^*\otimes x^*$. Then the involution on $\CBsAA(\haaga,A)\cong A$ is given by $(\psi^*)(u)=(\psi(u^*))^*$,
$u\in\haaga$, and similarly for $\CBsAA(M_I,A)$ noting that $M_I$ is *-invariant.
\end{remark}

\medskip

\enddocument